\documentclass[11pt]{amsart}

\usepackage{graphicx,calc}

\newtheorem{theorem}{Theorem}[section]
\newtheorem{corollary}[theorem]{Corollary}
\newtheorem{proposition}[theorem]{Proposition}
\newtheorem{lemma}[theorem]{Lemma}

\newcommand{\longpage}{\enlargethispage{\baselineskip}}
\newcommand{\shortpage}{\enlargethispage{-\baselineskip}}

\newcommand{\mapright}[1]{\stackrel{#1}{\longrightarrow}}

\newcommand{\Aut}{\operatorname{Aut}}
\newcommand{\I}{\operatorname{I}}
\newcommand{\E}{\operatorname{\mathcal{E}}}
\newcommand{\Homeo}{\operatorname{Homeo}}

\newcommand{\PSL}{\operatorname{PSL}}
\newcommand{\Sz}{\operatorname{Sz}}
\newcommand{\W}{\operatorname{\mathcal{W}}}
\newcommand{\Z}{\operatorname{\mathbb{Z}}}

\begin{document}

\title{Free actions on handlebodies}

\author{Darryl McCullough}
\address{Department of Mathematics\\
University of Oklahoma\\
Norman, Oklahoma 73019\\
USA}
\email{dmccullough@math.ou.edu}
\urladdr{www.math.ou.edu/$_{\widetilde{\phantom{n}}}$dmccullo/}

\author{Marcus Wanderley}
\address{Departmento de Matematica\\
Universidade Federal de Pernambuco\\
Av.\ Prof.\ Luiz Freire, s/n\\
Cid. Universitaria-Recife-PE\\
CEP 50.740-540, Brazil}
\email{mvw@dmat.ufpe.br}


\date{\today}


\begin{abstract}
The equivalence (or weak equivalence) classes of orien\-ta\-tion-preserving
free actions of a finite group $G$ on an orientable $3$-dimen\-sion\-al
handlebody of genus $g\geq 1$ can be enumerated in terms of sets of
generators of $G$. They correspond to the equivalence classes of generating
$n$-vectors of elements of $G$, where $n=1+(g-1)/|G|$, under Nielsen
equivalence (or weak Nielsen equivalence). For abelian and dihedral $G$,
this allows a complete determination of the equivalence and weak
equivalence classes of actions for all genera. Additional information is
obtained for solvable groups and for the groups $\PSL(2,3^p)$ with $p$
prime.  For all $G$, there is only one equivalence class of actions on the
genus $g$ handlebody if $g$ is at least $1+\ell(G)\,|G|$, where $\ell(G)$
is the maximal length of a chain of subgroups of $G$.  There is a
stabilization process that sends an equivalence class of actions to an
equivalence class of actions on a higher genus, and some results about its
effects are obtained.
\end{abstract}

\maketitle

\section*{Introduction}
\label{sec:intro}

The orientation-preserving free actions of a finite group $G$ on
3-dimension\-al orientable handlebodies have a close connection with a
long-studied concept from group theory, namely \emph{Nielsen equivalence}
of generating sets. Indeed, as we observe in section~\ref{sec:theory}
below, the free actions of $G$ on the handlebody of genus $g$, up to
equivalence, correspond to the Nielsen equivalence classes of $n$-element
generating sets of $G$, where $n=1+(g-1)/|G|$. We will utilize this to
prove a number of results about equivalence and weak equivalence of free
actions. These results are summarized in concise form in
section~\ref{sec:results}, which also contains precise definitions of
equivalence, weak equivalence, and other concepts that we shall use.

A special feature of free actions on handlebodies is that there is a
stabilization process relating actions on different genera. When a
handlebody $V$ with a free $G$-action contains a $G$-invariant handlebody
$U$ such that $\overline{V-U}$ consists of disjoint $1$-handles, the action
on $V$ is called a stabilization of the action on $U$. Inequivalent actions
can become equivalent after stabilization, indeed we do not know an example
of actions that remain inequivalent after even an elementary stabilization
(i.~e.~a stabilization for which $\overline{V-U}$ consists of $|G|$
$1$-handles). Such an example could not involve a solvable group, since a
result of M. Dunwoody implies that for solvable $G$, any two actions on a
handlebody of more than the minimum possible genus for a $G$-action are
equivalent (see corollary~\ref{coro:solvable corollary 2} below). For an
arbitrary $G$, proposition~\ref{prop:stabilization} shows that any two
actions become equivalent after $\mu(G)$ elementary stabilizations, where
$\mu(G)$ is the minimum number of generators of $G$. We remark that there
is an interesting theory of stabilization of actions on
$2$-manifolds~\cite{Kulkarni, M-M}.

The connection between free actions on handlebodies and Nielsen equivalence
is well-known in some circles, although we cannot find an explicit
statement in the literature.  It was known to J. Kalliongis and A. Miller
and is is a direct consequence of theorem~1.3 in their paper \cite{K-Mi1}
(for free actions, the graph of groups will have trivial vertex and edge
groups, and the equivalence of graphs of groups defined there is readily
seen to be the same as Nielsen equivalence on generating sets of $G$).
Indeed, more delicate classifications of nonfree actions on handlebodies
have been examined in considerable depth.  A general theory of actions was
given in~\cite{M-M-Z} and ~\cite{K-Mi1}, and the actions on very low genera
were extensively studied in~\cite{K-Mi2}. Actions with the genus small
relative to the order of the group are investigated in~\cite{M-Z}, and the
special case of orientation-reversing involutions is treated in~\cite{K-M}.
The first focus on free actions seems to be~\cite{Przytycki}, where it was
proven that for a cyclic group, any free action on a handlebody of genus
above the minimal one is the stabilization of an action on minimal genus,
and that any two free actions on a handlebody are weakly equivalent. These
results were generalized to dihedral and abelian groups in~\cite{W1,W2},
whose results are reconfirmed and extended in section~\ref{sec:abelian and
dihedral} below.

Some of the arguments in this paper can be shortened by invoking results
from~\cite{M-M-Z}. Since the general theory given there is much more
elaborate than the elementary methods needed for the present work, we have
chosen to make our arguments self-contained.  After giving some more
precise definitions and stating our main results in
section~\ref{sec:results}, we develop the general theory relating free
actions to Nielsen equivalence in section~\ref{sec:theory}. We apply this
in section~\ref{sec:solvable groups} to treat the case when $G$ is
solvable, and in section~\ref{sec:abelian and dihedral} we examine the
specific cases of abelian and dihedral groups.  In section~\ref{sec:Free
Actions of PSL(2,3p)}, we show that for $p$ prime, two free actions of
$\PSL(2,3^p)$ on a handlebody of genus above the minimal one are
equivalent.  By work of Evans \cite{Evans1} and Gilman \cite{Gilman}, it is
known that the same is true for $\PSL(2,2^m)$ (for all $m\geq 2$) and
$\PSL(2,p)$ (for $p$ prime).  Section~\ref{sec:stabilization of actions}
gives some general results on stabilization, in particular, we prove that
if the genus of $V$ is at least $1+\ell(G)|G|$, where $\ell(G)$ is the
maximum length of a decreasing chain of nonzero subgroups of $G$, then any
two free $G$-actions are equivalent. In section~\ref{sec:questions}, we
state some open problems. In particular, do there even exist inequivalent
actions that are not minimal genus actions?

\section[Results]{Statement of results}\label{sec:results}
Two (effective) actions $\rho_1,\rho_2\colon G\to\Homeo(X)$ are said to be
\emph{equivalent} if they are conjugate as representations, that is, if
there is a homeomorphism $h\colon X\to X$ such that
$h\rho_1(g)h^{-1}=\rho_2(g)$ for each $g\in G$. They are \emph{weakly
equivalent} if their images are conjugate, that is, if there is a
homeomorphism $h\colon X\to X$ so that $h\rho_1(G)h^{-1}=\rho_2(G)$. Said
differently, there is some automorphism $\alpha$ of $G$ so that
$h\rho_1(g)h^{-1} = \rho_2(\alpha(g))$ for all $g$. In words, equivalent
actions are the same after a change of coordinates on the space, while
weakly equivalent actions are the same after a change of coordinates on the
space and possibly a change of the group by automorphism. If $X$ is
homeomorphic to $Y$, then the sets of equivalence (or weak equivalence)
classes of actions on $X$ and on $Y$ can be put into correspondence using
any homeomorphism from $X$ to $Y$.

Henceforth the term \emph{action} will mean an orientation-preserving
\emph{free} action of a finite group on a $3$-dimensional orientable
handlebody of genus $g\geq 1$ (only the trivial group can act freely on the
handlebody of genus~$0$, the $3$-ball). One may work with either
piecewise-linear or smooth actions; we assume that one of these categories
has been chosen, and that all maps, isotopies, etc.~lie in the category.

For a finite group $G$ we denote by $\mu(G)$ the minimum number of
generators in any generating set for $G$, and by $\ell(G)$ the maximum
$\ell$ such that $G=G_1\supset G_2\supset \cdots \supset
G_{\ell}\supset\{1\}$ is a properly descending chain of subgroups of
$G$. When there is only one group $G$ under consideration, we often write
$\mu$ for~$\mu(G)$ and $\ell$ for~$\ell(G)$.

Fix a finite group $G$, and consider an action of $G$ on a
handlebody $V$. The quotient map $V\to V/G$ is a covering map, so
the action corresponds to an extension
\begin{equation}
1\longrightarrow\pi_1(V)\longrightarrow\pi_1(V/G)
\mapright{\phi}G\longrightarrow1\ .
\tag{$*$}
\end{equation}

\noindent A torsionfree finite extension of a finitely generated free group
is free (by~\cite{K-P-S} any finitely generated virtually free group is the
fundamental group of a graph of groups with finite vertex groups, and if
the group is torsionfree, the vertex groups must be trivial). So
$\pi_1(V/G)$ is free. Since $V$ is irreducible, so is $V/G$, and
theorem~5.2 of~\cite{Hempel} shows that $V/G$ is a handlebody.  From
equation ($*$), the genus of $V/G$ must be at least~$\mu$, so its Euler
characteristic is at most $1-\mu$. Therefore the Euler characteristic of
$V$ is at most $|G|\,(1-\mu)$, and the genus of $V$ is at least
$1+|G|\,(\mu-1)$.  On the other hand, if $W$ is a handlebody of genus
$\mu$, we may fix any surjective homomorphism from $\pi_1(W)$ to $G$, and
the covering space $V$ corresponding to its kernel has genus
$1+|G|\,(\mu-1)$ and admits an action of $G$. So $1+|G|\,(\mu-1)$ is the
minimal genus among the handlebodies that admit a $G$-action. We denote
this minimal genus by~$\Psi(G)$, and we call an action of $G$ on a
handlebody of genus $\Psi(G)$ a \emph{minimal genus action.}

If $V$ is any handlebody that admits a $G$-action, then the genus of $V/G$
is $\mu + k$ for some integer $k\geq 0$, so the genus of $V$ is
$1+|G|\,(\mu + k -1)$. Denote by $\E(k)$ the set of equivalence classes of
actions of $G$ on a handlebody of genus $1+|G|\,(\mu + k -1)$, and by
$\W(k)$ the set of weak equivalence classes. In particular, $\E(0)$ and
$\W(0)$ are the equivalence classes of minimal genus actions.  There is a
natural surjection from $\E(k)$ to~$\W(k)$. In the introduction we
explained that an action of $G$ on $V$ is a stabilization of an action on
$U$ if and only if there is a $G$-equivariant imbedding of $U$ into $V$
such that $\overline{V-U}$ consists of $1$-handles. We now give a more
convenient description. For $G$ acting on $V$, choose a disc
$E\subset\partial V$ which is disjoint from all of its
$G$-translates. Attach a $1$-handle $D^2\times \I$ to $V$ using an
orientation-reversing imbedding $j\colon D^2\times\partial \I\to E$. For
each $g\in G$ attach a $1$-handle to $g(E)$ using the imbedding $g\circ
j$. The $G$-action on $V$ extends to an action of $G$ on the union of $V$
with these $1$-handles, and this action is called an \emph{elementary
stabilization} of the original action. Alternatively, we may think of this
as attaching a $1$-handle to a disc in the boundary of $V/G$ and extending
the homomorphism $\phi\colon \pi_1(V/G)\to G$ that determines the action in
the exact sequence $(*)$ to a homomorphism from $\pi_1(V/G)*\Z$ to $G$ by
sending the generator of $\Z$ to~$1$. Since any two discs in $\partial V/G$
are isotopic, and any orientation-reversing imbeddings of
$D^2\times\partial \I$ into a disc are isotopic, the equivalence class of
the resulting action is well-defined. The result of applying some number of
elementary stabilizations is called a \emph{stabilization} of the original
action.  For each $k\geq 0$ and $m\geq1$, stabilization gives well-defined
functions from $\E(k)$ to $\E(k+m)$ and from $\W(k)$ to $\W(k+m)$.
\longpage

Let $e(k)$ denote the cardinality of $\E(k)$, and $w(k)$ the cardinality
of~$\W(k)$. Clearly~$e(k)\geq w(k)$, and $w(k)\geq 1$ for $k\geq 0$. Recall
that the Euler $\varphi$-function is defined on positive integers by
$\varphi(1)=1$, and for $m>1$, $\varphi(m)$ is the number of integers $q$
with $1\leq q<m$ such that $\gcd(m,q)=1$. Notice that for $m>2$,
$\varphi(m)$ is even, since if $\gcd(m,q)=1$ then $\gcd(m,m-q)=1$. Using
these notations, we can give concise statements of our main results:\par
\begin{enumerate}
\setlength{\itemsep}{1 ex}
\setlength{\parskip}{0.5 ex}
\item If $G$ is solvable, then $e(k)=1$ for all $k\geq 1$
(corollary~\ref{coro:solvable corollary 2}), while $w(0)$ can be
arbitrarily large (theorem~\ref{thm:N minimal classes}).
\item If $G$ is abelian, with $G \cong \Z/d_1\oplus\cdots \oplus \Z/d_m$,
where $d_{i+1}|d_i$ for $1\leq i<m$, then $w(0)=1$, and $e(0)=1$ if
$d_m=2$, otherwise $e(0)=\varphi(d_m)/2$ (theorem~\ref{thm:abelian}).
\item If $G$ is dihedral of order~$2m$, then $w(0)=1$, and $e(0)=1$ if
$m=2$, otherwise $e(0)=\varphi(m)/2$ (theorem~\ref{thm:dihedral}).
\item If $G=\PSL(2,3^p)$ with $p$ prime, then $e(k)=1$ for all $k\geq 1$
(corollary~\ref{coro:actions of PSL(2,3p)}).
\item The smallest genus of handlebody admitting inequivalent actions of a
group is $g=11$, which has two equivalence classes of actions of the
dihedral group of order $10$ (corollary~\ref{coro:smallest genus actions}).
\item The smallest genus of handlebody admitting inequivalent actions of an
abelian group is $g=26$, which has two equivalence classes of $C_5\times
C_5$-actions (corollary~\ref{coro:smallest genus actions}).
\item For all $G$ and all $k\geq 0$, $\E(k)\to \E(k+\mu(G))$ and $\W(k)\to
\W(k+\mu(G))$ are constant (proposition~\ref{prop:stabilization}).
\item For all $G$ and all $k> \ell(G)-\mu(G)$, $e(k)=1$
(corollary~\ref{coro:limit stabilization}).
\end{enumerate}

\noindent
In addition to these general results, we will see calculations for several
specific groups.

\begin{figure}
\includegraphics[width=\textwidth-0.2in]{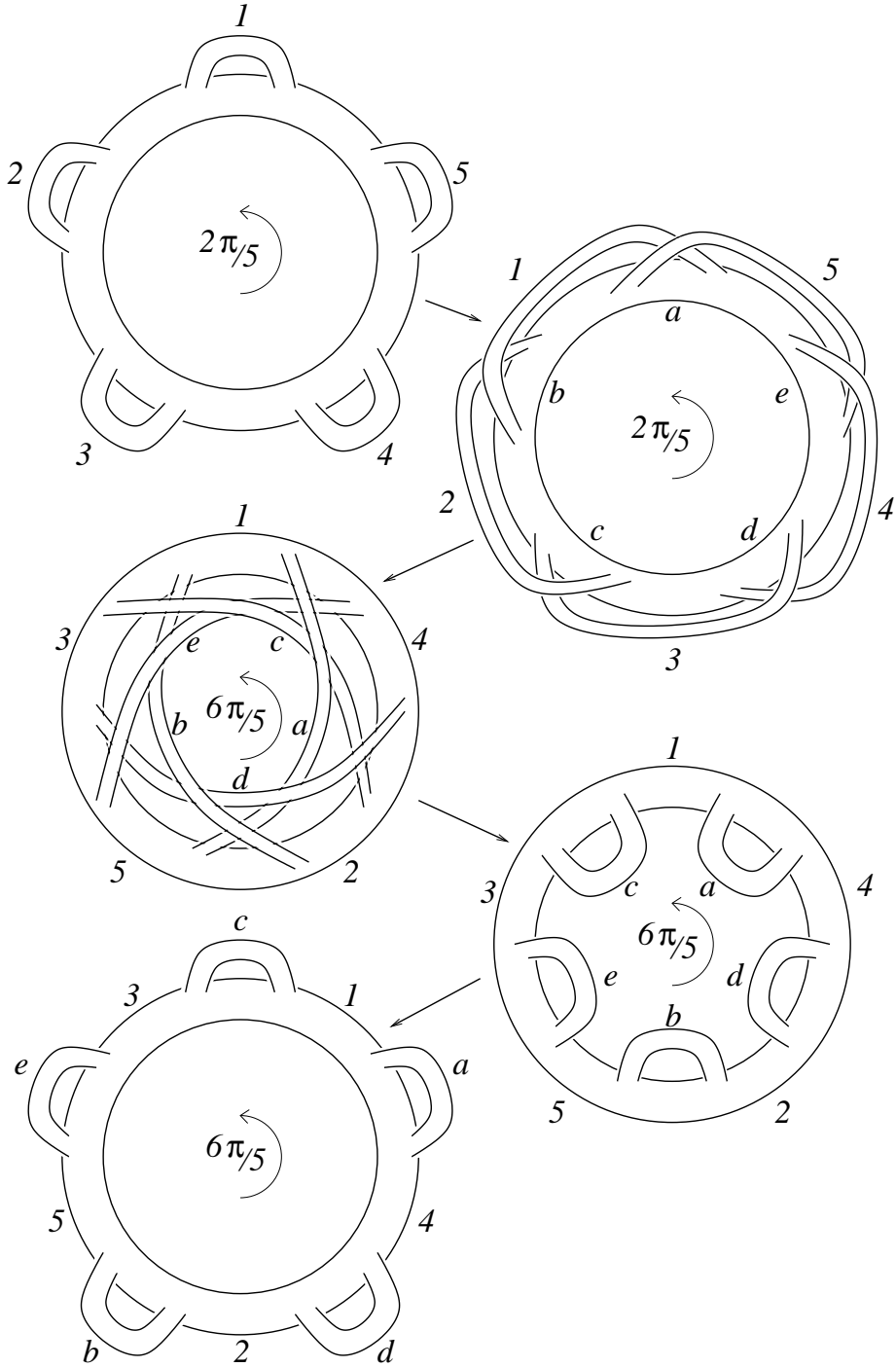}
\caption{An equivalence between the stabilizations of inequivalent actions
on the solid torus.}
\label{fig:actions}
\end{figure}

The results given above indicate that there is a strong tendency for
actions to become equivalent after stabilization, and as mentioned in the
introduction, we do not know even a single example of inequivalent actions
which do not become equivalent after one elementary stabilization. Even in
the simplest of cases, however, the underlying topology of an equivalence
between the stablizations of two inequivalent actions can be surprisingly
complicated. Figure~\ref{fig:actions} gives the steps of a visualization of
an explicit equivalence between two actions of the cyclic group $C_5$ of
order~$5$ on the handlebody of genus~$6$ which are the stabilizations of
inequivalent actions of $C_5$ on the solid torus (the classification of the
actions of cyclic groups is detailed in corollary~\ref{coro:cyclic
actions}).

\section[Nielsen equivalence]{Free actions and Nielsen equivalence}
\label{sec:theory}
The main result of this section, theorem~\ref{thm:equivalence classes},
gives the correspondence between actions of $G$ and Nielsen equivalence
classes of generating vectors of $G$. The proof is elementary, requiring
only the basic theory of covering spaces and some well-known facts about
free groups. We deduce, in corollary~\ref{coro:stabilization criterion}, an
algebraic criterion for an action to be a stabilization of another
action. We close with a few examples that illustrate the theory.

Consider two actions of $G$ on genus $g$ handlebodies $V_1$ and
$V_2$. Since the Euler characteristics of the quotients must be equal (to
$(1-g)/|G|$), we may choose diffeomorphisms from them to a single
handlebody $W$. The original actions are then equivalent to the actions of
$G$ by covering transformations on the covering spaces of $W$ determined by
two surjective homomorphisms $\phi_1,\phi_2\colon\pi_1(W)\to G$.  That is,
in classifying actions of a fixed group $G$ on handlebodies of a fixed
genus $g$, up to equivalence or up to weak equivalence, we may assume that
their quotients are the same handlebody~$W$.

A \emph{generating vector} for a group is a tuple $S=(s_1,\ldots,s_n)$ such
that $\{s_1,\ldots,s_n\}$ is a generating set.  A generating vector
$T=(t_1,\ldots,t_n)$ of $H$ is said to be obtained from $S$ by a
\emph{Nielsen move} if there is a $j$ so that $s_i=t_i$ for all $i\neq j$,
and for some $k\neq j$, $t_j$ equals $s_js_k^{\pm1}$ or
$s_k^{\pm1}s_j$. Also, an interchange of two of the entries, or the
replacement of an entry by its inverse, is a Nielsen move. Generating
vectors are called \emph{Nielsen equivalent} if there is a sequence of
Nielsen moves that changes one to the other. Generating vectors
$S=(s_1,\ldots,s_n)$ and $T=(t_1,\ldots,t_n)$ are called \emph{weakly
Nielsen equivalent} if there is an automorphism $\alpha$ of $G$ such that
$\alpha(S)$ is equivalent to $T$, where
$\alpha(S)=(\alpha(s_1),\ldots,\alpha(s_n))$. Note that equivalent or
weakly equivalent generating vectors must have the same number of
elements. A minimal generating vector for a free group is called a
\emph{basis.}

If $H$ is a free group with basis $S$, then any Nielsen move on $S$ induces
an automorphism of $H$. Nielsen proved \cite{Nielsen} that any two bases
for $H$ are Nielsen equivalent, consequently the Nielsen moves generate the
automorphism group of $H$.  Associated to a given handlebody structure on a
handlebody $W$ is a standard basis of the free group $\pi_1(W)$, where the
$i^{th}$ generator corresponds to a loop that goes once around the $i^{th}$
$1$-handle and not around any other handle. Any Nielsen move on this basis
is induced by an orientation-preserving diffeomorphism of $W$ (see for
example~\cite{M-M1}), so any automorphism of $\pi_1(W)$ can be induced by
an orientation-preserving diffeomorphism of~$W$. This also shows that
associated to any basis of $\pi_1(W)$ is a handlebody structure with
respect to which the basis is standard.

\begin{lemma} Let $W$ be a handlebody, and let $\phi_1,\phi_2\colon
\pi_1(W)\to G$ be surjective homomorphisms to a finite group $G$. Let
$(X_1,\ldots,X_n)$ be a basis for $\pi_1(W)$, so that
$S=(\phi_1(X_1),\ldots,\phi_1(X_n))$ and
$T=(\phi_2(X_1),\ldots,\phi_2(X_n))$ are generating vectors for $G$. Then
$S$ and $T$ are weakly Nielsen equivalent if and only there are an
isomorphism $\psi\colon \pi_1(W)\to\pi_1(W)$ and an isomorphism $\alpha$ of
$G$ such that $\alpha\phi_1=\phi_2\psi$.  They are Nielsen equivalent if
and only if $\alpha$ can be taken to be the identity automorphism of~$G$.
\label{lem:Nielsen equivalence}
\end{lemma}

\begin{proof} Suppose $\psi$ and $\alpha$ exist.  Since
$(\psi(X_1),\ldots,\psi(X_n))$ is a basis for $\pi_1(W)$, there is a
sequence of Nielsen moves that carries $(\psi(X_1),\ldots,\psi(X_n))$ to
$(X_1,\ldots,X_n)$. Applying $\phi_2$ shows that the corresponding Nielsen
moves in $G$ carry $(\phi_2\psi(X_1),\ldots,\phi_2\psi(X_n))$ to
$(\phi_2(X_1),\ldots,\phi_2(X_n))$. The latter is $T$, and since
$\phi_2\psi(X_i)=\alpha\phi_1(X_i)$, the former is $\alpha(S)$. Conversely,
suppose the generating vectors are weakly Nielsen equivalent. For a
sequence of Nielsen moves carrying $T$ to $\alpha(S)$, lifting the Nielsen
moves to corresponding Nielsen moves starting from $(X_1,\ldots,X_n)$
yields an isomorphism $\psi$ of $\pi_1(W)$ carrying $(X_1,\ldots,X_n)$ to a
basis $(\psi(X_1),\ldots,\psi(X_n))$ such that
$\phi_2(\psi(X_i))=\alpha\phi_1(X_i)$. This proves the lemma for weak
equivalence. The same argument, taking $\alpha$ to be the identity
automorphism, gives it for equivalence.\end{proof}

To translate this algebraic information into statements about equivalence
of actions, we use the following consequence of the theory of covering
spaces.
\begin{lemma} Let $\phi_1,\phi_2\colon \pi_1(W)\to G$
determine actions on handlebodies. The actions are weakly equivalent if and
only if there are an isomorphism $\psi\colon\pi_1(W)\to\pi_1(W)$ and an
isomorphism $\alpha\colon G\to G$ such that $\alpha\phi_1=\phi_2\psi$. They
are equivalent if and only if $\alpha$ may be taken to be the identity.
\label{lem:equivalence}
\end{lemma}

\begin{proof} Suppose that $\psi$ and $\alpha$ exist. For $j=1,2$, let
$V_j$ be the covering space of $W$ corresponding to the kernel of
$\phi_j$. Choose a diffeomorphism $h\colon W\to W$ inducing $\psi$. Since
$\psi$ must take the kernel of $\phi_1$ to the kernel of $\phi_2$, $h$
lifts to a diffeomorphism from $V_1$ to $V_2$. Also, $\psi$ induces
$\alpha$ from $G=\pi_1(W)/\hbox{ker}(\phi_1)$ to
$G=\pi_1(W)/\hbox{ker}(\phi_2)$. Using covering space theory, one can check
that this determines a weak equivalence of the actions, and an equivalence
when $\alpha$ is the identity automorphism.  Conversely, if the weak
equivalence (or equivalence, when $\alpha=1$) $H\colon V_1\to V_2$ exists,
then since $H(gx)=\alpha(g)H(x)$, $H$ induces a diffeomorphism $h\colon
W\to W$. Again by covering space theory, the induced automorphism $\psi$ of
$h$ on $\pi_1(W)$ satisfies $\phi_2\psi=\alpha\phi_1$.
\end{proof}

Putting these two lemmas together gives our main classification theorem.
\begin{theorem}
Let $G$ be finite, let $n\geq1$, and let $g= 1+|G|\,(n-1)$. Then the weak
equivalence classes of actions of $G$ on genus $g$ handlebodies correspond
bijectively to the weak equivalence classes of $n$-element generating sets
for $G$. The equivalence classes of actions correspond to the equivalence
classes of $n$-element generating sets.
\label{thm:equivalence classes}
\end{theorem}

\begin{proof} For $n<\mu$, both sets are empty, so we assume that $n\geq \mu$.
Let $W$ be a handlebody of genus $n$ and fix a basis
$(X_1,\ldots,\allowbreak X_n)$ for $\pi_1(W)$.  Any action is equivalent to
one determined by a surjective homomorphism $\phi\colon\pi_1(W)\to G$.  By
lemmas~\ref{lem:equivalence} and~\ref{lem:Nielsen equivalence}, sending the
action determined by $\phi$ to the basis $(\phi(X_1),\ldots,\phi(X_n))$
determines a bijection from the set of (equivalence or) weak equivalence
classes to the set of (equivalence or) weak equivalence classes of
$n$-element generating vectors of $G$.
\end{proof}

Let $\phi\colon \pi_1(W)\to G$ determine a $G$-action on a
handlebody $V$, and let $W'$ be obtained from $W$ by attaching a
single $1$-handle. Then $\pi_1(W')=\pi_1(W)*\Z$, where the
generator of $\Z$ corresponds to a loop that goes once around the
additional $1$-handle. We have noted that the action resulting
from an elementary stabilization of the action on $V$ is
determined by the homomorphism $\phi'\colon \pi_1(W')\to G$ which
equals $\phi$ on $\pi_1(W)$ and sends the generator of $\Z$
to~$1$.

\begin{corollary}
Suppose $G$ acts on handlebodies $U$ and $V$, where the genus of $U$ is
less than the genus of $V$. Let $(X_1,\ldots, X_m)$ be a basis for
$\pi_1(U/G)$ and let $(Y_1,\ldots,Y_n)$ be a basis for $\pi_1(V/G)$, and
let $\phi_U\colon \pi_1(U/G)\to G$ and $\phi_V\colon \pi_1(V/G)\to G$
determine the actions. Then the action on $V$ is a stabilization of the
action on $U$ if and only if the generating $n$-vectors
$(\phi_U(X_1),\ldots,\phi_U(X_m),1,\ldots,1)$ and
$(\phi_V(Y_1),\allowbreak\ldots,\allowbreak\phi_V(Y_n))$ of $G$ are
equivalent.
\label{coro:stabilization criterion}
\end{corollary}

\begin{proof}
Put $W_U=U/G$ and $W_V=V/G$. Since for any basis there is a handlebody
structure for which the basis is standard, there is an inclusion $j\colon
W_U\to W_V$ so that $\pi_1(W_U)\to \pi_1(W_V)$ carries $X_i$ to $Y_i$ for
$i\leq m$, and so that $\overline{W_V-j(W_U)}$ is a disjoint union of
$1$-handles. Then, the stabilized action is determined by the homomorphism
$\phi_V'\colon \pi_1(W_V)\to G$ for which $\phi_V'(Y_i)=\phi_U(X_i)$ for
$i\leq m$, and $\phi_V'(Y_i)=1$ for $i>m$. By theorem~\ref{thm:equivalence
classes}, the stabilized action and the action on $V$ are equivalent if and
only if $(\phi_U(X_1),\ldots,\phi_U(X_m),1,\ldots,1)$ and
$(\phi_V(Y_1),\allowbreak\ldots,\allowbreak\phi_V(Y_n))$ of $G$ are Nielsen
equivalent.
\end{proof}

This provides a simple criterion for an action to be a stabilization. A
vector of elements of $G$ is called \emph{redundant} if it is Nielsen
equivalent to a vector with an entry equal to~$1$. Notice that a vector is
redundant if it is even weakly Nielsen equivalent to a vector with an entry
equal to~$1$.

\begin{corollary} An action of $G$ on a handlebody $V$ is a
stabilization if and only if a generating vector $(s_1,\dots,s_n)$ of $G$
corresponding to the action as in theorem~\ref{thm:equivalence classes} is
redundant.  
\label{coro:redundant vector}
\end{corollary}

The literature contains a number of results on Nielsen equivalence of
generating vectors.  They were used to count $G$-defining subgroups of free
groups in \cite{Hall}. The action of the automorphism group of the free
group on generating vectors appears to have been introduced
in~\cite{Neumann-Neumann}. For infinite groups, there are quite a few
instances of inequivalent generating vectors of cardinality greater than
$\mu$. The paper of M. Evans~\cite{Evans2} gives general constructions of
these, as well as a summary of earlier results. But for finite groups, no
such example is known (see section~\ref{sec:questions} below). In the
remainder of this section, we will collect some of the known calculations
for specific finite groups, and give their consequences for group
actions. Some important general results of M.~Dunwoody for solvable groups
will be stated and used in the next section.

For $G=A_5$, B. Neumann and H. Neumann~\cite{Neumann-Neumann} (see
also~\cite{Stork1}) showed that there are two weak equivalence classes of
generating $2$-vectors for $A_5$.  Thus there are two weak equivalence
classes of actions on the handlebody of genus~$\Psi(A_5)=1+60(2-1)=61$.
D. Stork~\cite{Stork2} carried out similar calculations for
$\hbox{PSL}(2,7)$ and $A_6$. These show, for example, that there are $4$
weak equivalence classes of $A_6$-actions on the handlebody of genus $361$,
the minimal genus. Techniques developed by M. Lustig~\cite{Lustig} using
the Fox calculus yield additional examples of inequivalent generating
systems for certain groups.

Gilman~\cite{Gilman} proved that for $p$ prime, all $3$-element generating
sets of $\PSL(2,p)$ are equivalent. That is, $e(k)=1$ for all $k\geq 1$ for
these groups. In particular, this holds for $\PSL(2,5)\cong A_5$, so the
two equivalence classes of actions of $A_5$ on the handlebody of genus~$61$
become equivalent after a single elementary stabilization.
In~\cite{Evans1}, M. Evans proved that when $G$ is $\hbox{PSL}(2,2^m)$ or
the Suzuki group $\Sz(2^{2m-1})$ for $m\geq 2$, then for any $n\geq3$, all
$n$-element generating vectors are equivalent. Since $\mu=2$ for any of
these groups, this says that all actions of one of these groups on any
genus above the minimal genus are equivalent. In section~\ref{sec:Free
Actions of PSL(2,3p)}, we will prove a similar result for $\PSL(2,3^p)$
with $p$ prime. This includes the case of $\PSL(2,9)\cong A_6$, so the four
inequivalent $A_6$-actions on the handlebody of genus~$361$ all become
equivalent after a single elementary stabilization.

\section[Solvable Groups]{Actions of solvable groups}\label{sec:solvable
groups} 
In this section we use results of Dunwoody to examine the actions
of solvable groups. They show that although there can be an arbitrarily
large number of weak equivalence classes of minimal genus actions, all
actions on a handlebody whose genus is above the minimal one are
equivalent.

From~\cite{Dunwoody1}, we have the following fact.

\begin{theorem}
Let $G$ be a solvable group, and let $n>\mu(G)$. Then any two $n$-element
generating vectors are Nielsen equivalent.
\label{thm:Dunwoody theorem}
\end{theorem}

\noindent Applying theorem~\ref{thm:equivalence classes}, we have
immediately:

\begin{corollary}
Let $G$ be solvable and let $g$ be greater than $\Psi(G)$. Then any two
actions of $G$ on handlebodies of genus $g$ are equivalent.  Consequently,
any action of a solvable group on a handlebody not of minimal genus can be
destabilized to any given action of minimal genus.
\label{coro:solvable corollary 2}
\end{corollary}

Some examples due to Dunwoody also have implications for free actions.  The
following examples are from~\cite{Dunwoody2}:

\begin{theorem}
For every pair of positive integers $n$ and $N$, there exists a $p$-group
$G(n,N)$, nilpotent of length~$2$ and with $\mu(G(n,N))=n$, which has at
least $N$ weak Nielsen equivalence classes of $n$-element generating sets.
\label{thm:N minimal classes}
\end{theorem}

\noindent This yields immediately:

\begin{corollary}
For every pair of positive integers $n$ and $N$, there exists a $p$-group
$G(n,N)$, nilpotent of length~$2$, with $\mu(G(n,N))= n$ and $w(0)\geq
N$.
\label{coro:N minimal actions}
\end{corollary}

\noindent Of course, by corollary~\ref{coro:solvable corollary 2}, all of
the actions in $\E(0)$ become equivalent after a single stabilization.

We include here an elementary and transparent example of two weakly
inequivalent $G$-actions of a nilpotent group of order $2^{12}$ on the
handlebody of genus~$8193$. B. Neumann~\cite{Neumann} gave a nilpotent
group of order $2^{13}$ admitting weakly inequivalent actions on the
handlebody of this same genus; his example is slightly more complicated,
but requires only two generators, whereas ours requires~$3$. Dunwoody's
examples in~\cite{Dunwoody2} are considerably more sophisticated renderings
of the one we give here.

Let $G$ be the group with presentation
$$\langle x,y,z\;|\; x^8=y^8=z^{64}=1, [x,z]=[y,z]=1, [x,y]=z^8\rangle\ .$$

\noindent There is an extension$$1\to C_{64}\to G\to C_8\oplus C_8\to 1\ ,$$

\noindent where $C_{64}$ is the subgroup generated by $z$, and the
images of $x$ and $y$ generate $C_8\oplus C_8$.  Now
$xyx^{-1}=yz^8$, from which it follows that
$x^ay^bx^{-a}=y^bz^{8ab}$ for all integers $a$ and $b$.

We will show that the generating vectors $(x,y,z)$ and $(x,y,z^3)$
are not weakly Nielsen equivalent. By theorem~\ref{thm:equivalence
classes}, these correspond to two actions of $G$ on the handlebody
of genus $g=1+|G|\,(3-1)=8193$ which are not weakly equivalent.

Every element of $G$ can be written in the form $x^ay^bz^c$, where
$a$ and $b$ are integers mod 8, and $c$ is mod 64, and the inverse
of $x^ay^bz^c$ is $x^{-a}y^{-b}z^{-c-8ab}$. Using this, we can
calculate that $[x^py^qz^m,x^ry^sz^n]=z^{8(ps-rq)}$.

Sending $x$ to $(1,0,0)$, $y$ to $(0,1,0)$, and $z$ to $(0,0,1)$ defines a
homomorphism $G\to C_8\oplus C_8\oplus C_8$. Regarding these as vectors of
integers mod 8, any three-element generating set determines a $3\times 3$
matrix with entries mod 8. Nielsen moves on the generating set only change
the determinant of the associated matrix by multiplication by $\pm 1$. The
determinant of the matrix associated to $\{x,y,z^3\}$ is~$3$.

The subgroup $C_{64}$ is central. If $x^ay^bz^c$ has $0<a<8$, then
$y\,x^ay^bz^c\,y^{-1}=x^ay^bz^{c-8a}$, so $x^ay^bz^c$ is not
central. Similarly, if $0<b<8$, the element is not central. Therefore the
center of $G$ is exactly $C_{64}$, and any automorphism of $G$ must carry
$z$ to $z^d$ for some~$d$.

Now consider any automorphism $\alpha$ of $G$, with
$\alpha(x)=x^py^qz^m$, $\alpha(y)=x^ry^sz^n$, and $\alpha(z)=z^d$. From
above, $\alpha([x,y])=z^{8(ps-rq)}$.  Since this must equal
$\alpha(z^8)=z^{8d}$, it follows that $ps-rq$ is congruent to $d$ modulo
$8$. The matrix associated to the generating vector
$(\alpha(x),\alpha(y),\alpha(z))$ is
$$\begin{pmatrix} p&q&m\\
                  r&s&n\\
                  0&0&d\end{pmatrix}\ ,$$

\noindent which has determinant $(ps-rq)\,d$. This is congruent to
$d^2\pmod{8}$. Since the only squares modulo $8$ are $1$ and $4$, it follows that
$(\alpha(x),\alpha(y),\alpha(z))$ cannot be Nielsen equivalent to the
generating set $(x,y,z^3)$. Therefore the generating vectors $(x,y,z)$ and
$(x,y,z^3)$ are not weakly Nielsen equivalent.

\section[Abelian and Dihedral Groups]{Abelian and dihedral groups}
\label{sec:abelian and dihedral}

In this section we examine the actions of abelian and dihedral groups. We
will see that for either of these two kinds of groups, all actions on the
minimal genus are weakly equivalent, but that there can be arbitrarily
large numbers of equivalence classes.

\begin{theorem}
Let $A$ be a finite abelian group, $G\cong \Z/d_1\oplus\cdots \oplus
\Z/d_m$ where $d_{i+1}|d_i$ for $1\leq i<m$.  Then $\Psi(A)=1+|A|(m-1)$. If
$g>\Psi(A)$, then any two $A$-actions on a handlebody of genus $g$ are
equivalent.  Any two $A$-actions on a handlebody of genus $\Psi(A)$ are
weakly equivalent. If $d_m=2$, then all $A$-actions on the handlebody of
genus $\Psi(A)$ are equivalent, while if $d_m>2$, then there are exactly
$\varphi(d_m)/2$ equivalence classes, where $\varphi$ denotes the Euler
$\varphi$-function.
\label{thm:abelian}
\end{theorem}

In preparation for the proof, we will prove a general lemma about
vectors of elements in solvable groups.  By a \emph{cyclic tower}
for a group $G$ we mean a descending sequence of subgroups
$G=G_1\supset G_2\supset \cdots\supset G_{m}\supset \{1\}$ such
that $G_{i+1}$ is normal in $G_i$ and $G_i/G_{i+1}$ is cyclic.  We
allow $G_i$ to equal $G_{i+1}$ for some $i$, also we define
$G_j=\{1\}$ for all~$j>m$.

\begin{lemma} Let $G$ be solvable and let
$G=G_1\supset G_2\supset \cdots \supset G_m\supset \{1\}$ be a
cyclic tower for $G$. Let $T=(t_1,\ldots,t_n)$ be a vector of
elements of $G$ (not necessarily a generating vector). Then $T$ is
Nielsen equivalent to a vector with $t_i\in G_i$ for all $i$ (in
particular, $t_i=1$ for any $i>m$).
\label{lem:cyclic towers}
\end{lemma}

\begin{proof} Let $H$ be the subgroup of $G$ generated by $T$. Replacing
each $G_i$ by $H\cap G_i$, we may assume that $T$ is a generating set
for~$G$.

Suppose first that $G$ is a cyclic group $C_m$, generated by $t$, and that
at least two of the $t_i$, say $t_1$ and $t_2$, are not equal to~$1$. Write
$t_1=t^a$ and $t_2=t^b$, where $1\leq a,b\leq m-1$. If $a\geq b$, replace
$t_1$ by $t_1t_2^{-1}$, then $t_1=t^{a-b}$ and $t_2=t^b$.  If $a<b$,
replace $t_2$ by $t_2t_1^{-1}$. Repeat this process until either $t_1=1$ or
$t_2=1$. By an interchange, we may assume that $t_1\neq 1$ and
$t_2=1$. Repeating the process with the other elements, we eventually
achieve that $t_i=1$ for all~$i\geq 2$.

In the general case, we may regard $T$ as a vector of elements in the
cyclic group $G_1/G_2$. By the cyclic case, we may assume that $t_i\in G_2$
for all $i>1$. The subgroup of $G$ generated by $\{t_2,\ldots,t_n\}$ has a
cyclic tower of length $m-1$. By induction on $m$, $(t_2,\ldots,t_n)$ is
Nielsen equivalent to a vector with~$t_i\in G_i$.
\end{proof}

\begin{proof}[Proof of Theorem~\ref{thm:abelian}] We have $\mu=m$,
since the minimal number of generators of $A\otimes (\Z/d_m)=(\Z/d_m)^m$ is
$m$, so $\Psi(A)=1+|A|(m-1)$. If $g>\Psi(A)$, then
corollary~\ref{coro:solvable corollary 2} shows that any two $A$-actions on
a handlebody of genus $g$ are equivalent.

Regard $G$ as solvable with cyclic tower given by
$G_i=\Z/d_i\oplus\cdots\oplus \Z/d_m$. Let $S=(s_1,\ldots,s_n)$ be any
generating vector. By lemma~\ref{lem:cyclic towers}, we may assume that
$s_i\in G_i$, and consequently $s_1$ generates $G_1/G_2$.  Now
$(s_2,\ldots,s_n)$ generate $G_2$, since otherwise the quotient of $A$ by
the subgroup that it generates is of the form $\Z/d_1\oplus A_2$ with $d_1$
and the order of $A_1$ not relatively prime, but this quotient could not be
generated by $s_1$. So we may apply Nielsen moves changing $s_1$ by
multiples of the other $s_i$ until $s_1\in \Z/d_1$. Inductively, we may
assume that $s_i\in \Z/d_i$ and generates $\Z/d_i$ for $i\leq m$, and
$s_i=1$ for $i>m$.

Let $T=(t_1,\ldots,t_m)$ be any another $n$-element generating vector for
$A$. Again, we may assume that $t_i\in \Z/d_i$ and $t_i$ generates
$\Z/d_i$, so $t_i=s_i^{p_i}$ and $s_i=t_i^{q_i}$ for some $p_i$ and $q_i$
relatively prime to $d_i$. By Nielsen moves, replace $t_2$ by
$t_2t_1^{q_1}$, then $t_1$ by
$t_1(t_2t_1^{q_1})^{-p_1}=t_1t_2^{-p_1}t_1^{-p_1q_1}=t_2^{-p_1}$.  Since
$p_1$ is relatively prime to $d_1$ and $d_2$ divides $d_1$, $p_1$ is also
relatively prime to $d_2$. So there exists $r$ with $(t_2^{-p_1})^r=t_2$,
and by Nielsen moves we may replace $t_2t_1^{q_1}$ by
$t_1^{q_1}=s_1$. Interchanging $t_1$ and $t_2$, we have that $t_1=s_1$ and
$t_2$ still generates $\Z/d_2$.  Continuing, we may assume that $t_i=s_i$
for all $i<n$. If $n>m$, then $t_i=s_i$ for all $i$ since both equal $1$
for $i>m$. This proves that all actions on genera greater than $\Psi(A)$
are equivalent. If $n=m$, then we have only that $t_m=s_m^p$, with $p$
relatively prime to~$d_m$.

If $d_m=2$, then $T$ must be Nielsen equivalent to $S$. From now on, assume
that $d_m>2$. The automorphism $\alpha$ of $A$ defined by $\alpha(s_i)=s_i$
for $i<m$ and $\alpha(s_m)=s_m^p$ shows that all $m$-element generating
vectors are weakly Nielsen equivalent. Since $(s_1,\ldots,s_{m-1},s_m^p)$
is Nielsen equivalent to $(s_1,\ldots,s_{m-1},s_m^{-p})$, there are at most
$\phi(d_m)/2$ equivalence classes. To show that this is a lower bound for
the number of equivalence classes, we define a function from (ordered)
generating sets to $m\times m$ matrices as follows. Regard $s_i$ as an
$m$-tuple with all entries $0$ except for a $1$ in the $i^{th}$
place. Working $\bmod{d_m}$, any generating $m$-vector then determines an
$m\times m$ matrix with $i^{th}$ row the vector corresponding to $t_i$. A
Nielsen move corresponds to multiplying by an elementary matrix (that adds
one row to another, or interchanges two rows, or multiplies one row by
$-1$). These elementary matrices have determinant $\pm1$, so Nielsen
equivalent generating sets have determinants that are either equal or are
negatives, as elements of $\Z/d_m$. The determinant of the matrix
corresponding to $(s_1,\ldots,s_{m-1},s_m^p)$ is $p$. This gives the lower
bound of $\phi(d_m)/2$ on the number of equivalence classes.
\end{proof}

Specializing to the case of a cyclic group, we have the following:

\begin{corollary}
For the cyclic group $C_k$, $\Psi(C_k)=1$. If $g>1$, then any two
$C_k$-actions on a handlebody of genus $g$ are equivalent.  On the solid
torus, any two $C_k$-actions are weakly equivalent, any two $C_2$-actions
are equivalent, and if $k>2$, then there are exactly $\varphi(k)/2$
equivalence classes of $C_k$-actions. 

Explicitly, if $\phi_q$ is the action on $S^1\times D^2$ defined by
$\phi_q(t)(\exp(i\theta),x)=(\exp(i\theta+2\pi iq/k),x)$ where $t$ is a
fixed generator of $C_k$ and $q$ is relatively prime to $k$, then
$\phi_{q_1}$ and $\phi_{q_2}$ are equivalent if and only if $q_1\equiv \pm
q_2 \pmod{k}$.\par
\label{coro:cyclic actions}
\end{corollary}

\begin{proof}
The corollary is immediate from the statement of theorem~\ref{thm:abelian},
except for the explicit description of the equivalence classes. The matrix
corresponding to $\phi_q$ in the proof of theorem~\ref{thm:abelian} is
$[q]$, and the last paragraph of the proof shows the condition for
equivalence of $\phi_{q_1}$ and~$\phi_{q_2}$.
\end{proof}

\begin{theorem}
Let $D_{2m}$ be the dihedral group of order $2m$. Then
$\Psi(D_{2m})=2m+1$. Any two $D_{2m}$-actions on the handlebody of genus
$2m+1$ are weakly equivalent. If $m=2$, all actions on the handlebody of
genus $2m+1$ are equivalent, and if $m>2$ then there are exactly
$\phi(m)/2$ equivalence classes.
\label{thm:dihedral}
\end{theorem}

\begin{proof} Regard $D_{2m}$ as
$\langle a,b\;|\; a^2=b^m=1,aba^{-1}=b^{-1}\rangle$, and let $C_m$ be the
cyclic subgroup of $D_{2m}$ generated by $b$. By lemma~\ref{lem:cyclic
towers}, applied to the tower $D_{2m}\supset C_m\supset\{1\}$, any
two-element generating vector $S=(x,y)$ is Nielsen equivalent to one of the
form $(ab^i,b^j)$. Since $ab^i$ has order~2, $b^j$ must generate $C_m$, so
$\gcd(m,j)=1$ and the vector is equivalent to $(a,b^j)$. For $m=2$, the proof
is complete. For $m>2$, there are at most $\phi(m)/2$ equivalence classes,
since $(a,b^j)$ is equivalent to $(a,b^{-j})$. On the other hand, it is
essentially an observation of D. Higman (see~\cite{Neumann}) that the pair
(possibly equal) of conjugacy classes of $[x,y]$ and $[y,x]$ is an
invariant of the Nielsen equivalence class of $(x,y)$. In our case,
$[a,b^j]^{\pm1}=b^{\mp 2j}$ Since $b^{2j}$ is conjugate to $b^{\pm2\ell}$
only when $j=\pm\ell$, there are exactly $\phi(m)/2$ equivalence classes.
\end{proof}

\begin{corollary} The smallest genus of handlebody admitting two
inequivalent actions is genus $11$, which has two equivalence classes of
$D_{10}$-actions. The smallest genus of handlebody admitting two
inequivalent actions of an abelian group is genus $26$, which has two
equivalence classes of $C_5\times C_5$-actions.
\label{coro:smallest genus actions}
\end{corollary}

\begin{proof}
Theorem~\ref{thm:abelian} verifies the assertion about abelian groups. By
theorem~\ref{thm:dihedral}, the smallest-genus inequivalent actions of
dihedral groups are the two $D_{10}$-actions in the corollary. Only cyclic
groups have $\mu(G)=1$, so $\Psi(G)\geq |G|+1$ for a noncyclic group. The
only nonabelian and nondihedral group of order smaller than $11$ is the
quaternion group, which is easily checked to have only one equivalence
class of generating pair.
\end{proof}

\section[Free Actions of $\PSL(2,3^p)$]{Free Actions of $\PSL(2,3^p)$ 
($p$ a prime number)}
\label{sec:Free Actions of PSL(2,3p)}

We have mentioned that Gilman~\cite{Gilman} proved that for $p$ prime, all
$3$-element generating sets of $\PSL(2,p)$ are equivalent, and
Evans~\cite{Evans1} proved the same for $\PSL(2,2^m)$ and for the Suzuki
groups $\Sz(2^{2m-1})$ for $m\geq 2$. In this section, we prove the same
for all $\PSL(2,3^p)$ with $p$ prime. The main result is the following.

\begin{theorem}
Let $p$ be prime. If $n>2$, then any $n$-element generating vector for
$\PSL(2,3^p)$ is redundant.
\label{thm:PSL(2,3p)}
\end{theorem}

\noindent
In particular, these groups include the case of $\PSL(2,9)\cong A_6$
(\cite{Suzuki}, p.~412). Before proving theorem \ref{thm:PSL(2,3p)}, we
deduce a corollary.
\begin{corollary}
Let $p$ be prime. If $n>2$, then any two $n$-element generating vectors for
$\PSL(2,3^p)$ are Nielsen equivalent. Consequently, for any handlebody of
genus above the minimal one, $1+3^p(3^{2p}-1)/2$, all $\PSL(2,3^p)$-actions
are equivalent.  
\label{coro:actions of PSL(2,3p)}
\end{corollary}

\begin{proof}
We recall from \cite{B-W} that a $2$-generator group $G$ is \emph{of
spread~$2$} when for any pair ${h_1,h_2}$ of nontrivial elements of $G$,
there exists $y\in G$ such that $\langle y,h_1\rangle=\langle
y,h_2\rangle=G$.  When $p>2$, theorem~4.02 of~\cite{B-W} shows that
$\PSL(2,3^p)$ has spread~$2$. For $p=2$, $\PSL(2,3^2)$ is isomorphic to
$A_6$ (\cite{Suzuki}, p.~412), so has spread~$2$ by proposition~3.02
of~\cite{B-W}.

Let $(s_1,\dots,s_n)$ and $(t_1,\dots,t_n)$ be any two $n$-element
generating sets. By repeated use of theorem~\ref{thm:PSL(2,3p)}, we may
assume they are of the form $(s_1,s_2,1,(1))$ and $(t_1,t_2,1,(1))$, where
$(1)$ indicates a possibly empty string of $1$'s. Choose $y$ such that
$\langle y,s_1\rangle=\langle y,t_1\rangle=\PSL(2,3^p)$. As in lemma~2.8
of~\cite{Evans1}, we have equivalences $(s_1,s_2,1,(1))\sim
(s_1,s_2,y,(1))\sim (s_1,y,1,(1))\sim (s_1,y,t_1,(1))\sim (y,t_1,1,(1))\sim
(y,t_1,t_2,(1)) \allowbreak\sim (t_1,t_2,1,(1))$.
\end{proof}

To prepare for the proof of theorem~\ref{thm:PSL(2,3p)}, we will list
several group-theoretic results. The first is lemma~4.10 of~\cite{Evans1}.

\begin{lemma}
Let $G$ be a simple group generated by involutions $g_1,\dots,g_n$, $n\geq
3$. Then $(g_1,\dots,g_n)$ is redundant.  
\label{lem:involutions}
\end{lemma}

By direct calculation, we have the following information.

\begin{lemma} Let $x$ and $y$ be elements of $S_4$ such that $x\not=y$ 
and $x\not=y^{-1}$.
\begin{enumerate}
\item[{\rm(a)}] If $|x|=|y|=3$, then $|xy|=|yx|=2+t$ and
$|xy^2|=|y^2x|=3-t$, where $t$ is either~$0$ or~$1$.
\item[{\rm(b)}] If $|x|=3$ and $|y|=4$, then $|xy|=|yx|=2+t$ and
$|xy^3|=|y^3x|=4-t$ where $t$ is either $0$ or $2$. Moreover,
$|xy^2|=|y^2x|=3$ and $|x^2y|=|yx^2|=2$.
\item[{\rm(c)}] If $|x|=|y|=4$, then $|xy|=3$ and
$|xy^2|=2$.
\end{enumerate}
\label{lem:equations}
\end{lemma}

\noindent Since $A_5\cong\PSL(2,5)$, theorem~$1$ of~\cite{Gilman} implies

\begin{lemma}
Let $T=(u,v,w)$ be a generating vector for $A_5$. Then $T$ is redundant.
\label{lem:A5 redundant}
\end{lemma}

\noindent Another useful property of $A_5$ follows
from~\cite{Neumann-Neumann}.

\begin{lemma}
Let $(u,v)$ be a generating vector for $A_5$, and let
$(a,b)\in\{(2,3),(2,5),(3,5),(5,5)\}$. Then $(u,v)$ is Nielsen equivalent
to a generating vector $(u',v')$ with $(|u'|,|v'|)=(a,b)$.
\label{lem:NN}
\end{lemma}

\begin{proof}
Let $F_2$ denote the free group on two letters $u$ and $v$.  According
to~\cite{Neumann-Neumann}, there is an automorphism $f\colon F_2\to F_2$
such that as elements of $A_5$, $(|f(u)|,|f(v)|)=(a,b)$. Since
$f\in\Aut(F_2)$, it can be written as a product of Nielsen moves on the
basis $(u,v)$, so as vectors of elements of $A_5$, $(u,v)$ is Nielsen
equivalent to $(f(u),f(v))$.
\end{proof}

We will use a nice observation from~\cite{Gilman}.
\begin{lemma}
Let $T=(u_1,u_2,\dots,u_n)$ be a generating vector for a nonabelian simple
group $G$, with $u_1,u_2\neq 1$. Then $T$ is equivalent to a generating
vector $(u_1',u_2,\dots,u_n)$ with $u_1'$ conjugate to $u_1$ and
$[u_1',u_2]\neq 1$.  
\label{lem:nonabelian}
\end{lemma}

\begin{proof}
The subgroup generated by the conjugates of $u_1$ by elements of $H=\langle
u_2,\dots,u_n\rangle$ is normal, hence equals $G$. Since $u_2$ cannot be
central, there must be some $w\in H$ such that $u_2$ does not commute with
$wu_1w^{-1}$, and $T$ is equivalent to $(wu_1w^{-1},u_2,\dots,u_n)$.
\end{proof}

We will use some well-known information about subgroups
of~$\PSL(2,q)$. Specializing theorems~3(6.25) and~3(6.26) of
Suzuki~\cite{Suzuki} to the case at hand, we obtain the following
description.

\begin{theorem}
Let $p$ be a prime. Then every subgroup of $\PSL(2,3^p)$ is contained in
one of the following kinds of subgroups:
\begin{enumerate}
\item[(a)] The dihedral groups of orders $3^p+1$ and $3^p-1$.
\item[(b)] A group $K$ of order $3^p(3^p-1)/2$.  A Sylow $3$-subgroup $Q$
of $K$ is elementary abelian, normal, and the factor group $K/Q$ is a
cyclic group of order $(3^p-1)/2$.
\item[(c)] $A_4$.
\item[(d)] $S_4$ and $A_5$, which are subgroups only when $p=2$.
\end{enumerate}
\label{thm:subgroups}
\end{theorem}

\noindent Notice that a Sylow $3$-subgroup $Q$ in
theorem~\ref{thm:subgroups} is also a Sylow $3$-subgroup of $\PSL(2,3^p)$,
since $|\PSL(2,3^p)|=3^p(3^{2p}-1)/2$, and the group $K$ is its normalizer,
since $K$ is maximal. We will refer to a subgroup $K$ as in
theorem~\ref{thm:subgroups}(b) as a \emph{Sylow $3$-normalizer.} 

We will also need the following facts about centralizers, which follow from
3(6.5) and 3(6.8) of~\cite{Suzuki}.

\begin{lemma}
The centralizer of any element in $\PSL(2,3^q)$ is a maximal abelian
subgroup, and is either cyclic of order prime to~$3$, isomorphic to
$C_2\times C_2$, or is a Sylow $3$-subgroup of $\PSL(2,3^q)$.  Distinct
maximal abelian subgroups have trivial intersection.  
\label{lem:centralizers}
\end{lemma}

\noindent In particular, any two distinct Sylow $3$-subgroups have trivial
intersection.

\begin{proof}[Proof of theorem~\ref{thm:PSL(2,3p)}] Fix an $n$-element 
generating vector $T=(t_1,\ldots,t_n)$ of $\PSL(2,3^p)$. Suppose first that
$n>3$. We will show that $T$ is redundant.

Put $H=\langle t_1,t_2,t_3\rangle$. If $H=\PSL(2,3^p)$, then $T$ is
redundant. If $H$ is cyclic, dihedral, $A_4$ or $S_4$, then $T$ is
redundant by theorem~\ref{thm:Dunwoody theorem}, and if $H$ is $A_5$, then
lemma~\ref{lem:A5 redundant} applies. If $H$ is contained in a Sylow
$3$-normalizer $K$ of a Sylow $3$-subgroup Q, then $(t_1,t_2,t_3)$ is
(Nielsen) equivalent to $(s_1,s_2,s_3)$ such that $s_1\in Q$ (and $s_1\neq
1$, otherwise $T$ is redundant). As $\PSL(2,3^p)\neq K$, there is a $j\geq
4$ for which $\langle s_1,s_2,t_j\rangle$ is not a subgroup of $K$. Since
the Sylow $3$-subgroups intersect trivially, $\langle s_1,s_2,t_j\rangle$
cannot be contained in any other Sylow $3$-normalizer.  So $\langle
s_1,s_2,t_j\rangle$ is redundant or $\PSL(2,3^p)=\langle
s_1,s_2,t_j\rangle$. In either case, $T$ is redundant.

In the rest of the proof, we may assume that $T=(u,v,w)$. We assume at
every stage of the argument that no two elements of $T$ generate a cyclic
group or all of $\PSL(2,3^p)$, since otherwise $T$ is clearly redundant.
We will argue that $T$ either is redundant or is equivalent to $(x,y,z)$
with $x$, $y$, and $z$ all of order~$2$. This will prove the theorem, since
lemma~\ref{lem:involutions} shows that the latter is also redundant.

We first show that $T$ is (redundant or) equivalent to $(x,v,w)$ with $x$
of order $2$. If none of $u$, $v$, or $w$ already has order~$2$, put
$H=\langle u,v\rangle$.  We may assume that $H$ is not contained in any
Sylow $3$-normalizer $K$. For if it is, then since $K/Q$ is cyclic, we may
assume that $u\in Q$ (and $u\neq 1$). Interchanging $v$ and $w$ if
necessary, we may assume that $H$ is not contained in $K$. Since the Sylow
$3$-subgroups intersect trivially, $H$ cannot be contained in any other
Sylow $3$-normalizer.

Since neither of its generators has order~$2$, $H$ cannot be dihedral. If
$H$ is $A_4$ or $S_4$, then lemma~\ref{lem:equations} can be used to find
$x$. If it is $A_5$, then lemma~\ref{lem:NN} can be used. So we may write
$T$ as $(x,v,w)$ with $x$ of order~$2$.

Next, we will show that $(x,v,w)$ is equivalent to $(x,y,w)$ with $y$ also
of order $2$. Put $H=\langle v,w\rangle$. By lemma~\ref{lem:nonabelian}, we
may assume that $H$ is nonabelian.

We may assume that $H$ is not contained in any Sylow $3$-normalizer. For
suppose it lies in $K$. Since $K/Q$ is cyclic, we may asssume that $v\in
Q$. Since $x$ cannot be in $K$, $xw$ cannot be in $K$. Since the Sylow
$3$-subgroups are disjoint, $\langle v,xw\rangle$ cannot be in any Sylow
$3$-normalizer. If $H$ is dihedral, then one of $v$ or $w$ already has
order~$2$ and will be $y$. If $H$ is one of $A_4$, $S_4$, or $A_5$, then as
before, either lemma~\ref{lem:equations} or~\ref{lem:NN} produces~$y$.

We now have $T=(x,y,w)$ where $x$ and $y$ are elements of order $2$.  We
may assume that $|w|>2$ and, using lemma~\ref{lem:nonabelian}, that
$[x,y]\neq 1$, so $|xy|>2$. Put $H=\langle xy,w\rangle$. It cannot be
dihedral, since neither generator is of order~$2$.

Suppose that $H$ is contained in a Sylow $3$-normalizer $K$.
If $xy$ is of order $3$, then since $x$ inverts $xy$, $x\in K$ and hence
$y\in K$, so $K=\PSL(2,3^p)$, a contradiction. So $xy\notin Q$.

Assume first that $p>2$, so that $|K|$ is odd. Then $\langle x,w\rangle$ is
not contained in a Sylow $3$-normalizer. If it is dihedral, then we can
take $z=xw$. Otherwise, it is isomorphic to $A_4$, so $w$ has
order~$3$. Similarly, considering $\langle x,(xy)w\rangle$, we may assume
that $(xy)w$ has order~$3$. Since this element lies in $K$, it must lie in
$Q$, forcing the contradiction that $xy\in Q$.

When $p=2$, we have from p.~398 of~\cite{Suzuki} that $K$ is a split
extension of $C_3\times C_3$ by $C_4$. If $k$ generates $C_4$, then $q^2$
must act by inverting each element of $C_3\times C_3$ (otherwise it would
act trivially, but then the centralizer of an element of $Q$ would be
larger than $Q$, contradicting lemma~\ref{lem:centralizers}) so $k^2q$ has
order~$2$ for any $q\in C_3\times C_3$.  Since $xy\notin Q$ and $|xy|>2$,
we have $|xy|=4$. If $w$ has order~$3$, then $w(xy)^2$ is of order $2$.  If
$w$ is of order~$4$, then $w(xy)$ is of order~$2$.  So we may assume that
$H$ is not contained in a Sylow $3$-normalizer.

If $H$ is isomorphic to $A_4$ or $S_4$, then lemma~\ref{lem:equations}
applies. The case of $p>2$ is complete, so we may assume that $p=2$ and
$H\cong A_5$. Since $T$ is equivalent to $(x,xy,w)$, we may apply
lemma~\ref{lem:NN} to $\langle xy,w\rangle$ to obtain a new generating
vector of the form $(x,y,w)$ where $x$ and $y$ have order $2$, $w$ has
order $5$, and $\langle y,w\rangle\cong A_5$. By the previous arguments, we
may also assume that $\langle x,w\rangle$ and $\langle xy,w\rangle$ are
isomorphic to~$A_5$.

In the remainder of the proof, it is convenient to regard $\PSL(2,9)$ as
$A_6$. Also, we may apply any automorphism to all elements of the
generating vector, since any vector weakly equivalent to a redundant vector
is redundant.

According to 3(2.19) of~\cite{Suzuki}, there are two conjugacy classes of
$A_5$-sub\-groups in $A_6$. One class consists of the six $A_5$-subgroups
that stabilize one of the six letters, and the other class consists of six
$A_5$-subgroups that act transitively on the six letters. Moreover, there
exists an automorphism $\alpha$ of $A_6$ that interchanges the two
classes. Applying $\alpha$ to each of the generating elements, if
necessary, we may assume that $\langle xy,w\rangle$ fixes a letter. Since
$w$ is a $5$-cycle, this must be the unique letter fixed by $w$. Neither of
$\langle x,w\rangle$ and $\langle y,w\rangle$ can fix a letter. For if so,
it would fix the same letter as $\langle xy,w\rangle$, and one of $x$ or
$y$ would be in $\langle xy,w\rangle$, forcing the latter to be all of
$A_6$. Applying $\alpha$ to all generators, $\langle x,w\rangle$ and
$\langle y,w\rangle$ become $A_5$-subgroups fixing the letter fixed by $w$,
so are equal, achieving the same contradiction.
\end{proof}
\shortpage

The principal difficulties in extending the proof of
theorem~\ref{thm:PSL(2,3p)} to more general cases $\PSL(2,q^r)$ seem to be
the more complicated Sylow $q$-normal\-izers, when $q\neq 3$, the presence of
subgroups of the form $\PSL(2,q^s)$ when $r$ is not prime, and analyzing
the case when any two generators generate an $A_5$ subgroup.

\section[Stabilization]{Stabilization of actions}\label{sec:stabilization
of actions} Throughout this section, $G$ is an arbitrary finite group, and
as usual $\mu(G)$ or just $\mu$ will denote the minimal number of
generators for $G$, and $\ell(G)$ or just $\ell$ the maximum length of a
chain of strictly decreasing nonzero subgroups of~$G$. Clearly
$\mu\leq\ell$, and if $|G|=p_1^{\alpha_1}\cdots p_r^{\alpha_r}$, then
$\ell\leq\alpha_1+\cdots +\alpha_r$.

\begin{proposition}
Any two actions of a group $G$ on handlebodies of the same genus become
equivalent after at most $\mu(G)$ stabilizations.
\label{prop:stabilization}
\end{proposition}

\begin{proof} Let $S=(s_1,\ldots,s_\mu)$ be a generating vector of
minimal length, and let $T=(t_1,\ldots,t_m)$ be any generating vector.  Put
$T'=(t_1,\ldots,t_m,1,\ldots,\allowbreak 1)$, where $\mu$ $1$'s have been
added. Then $T'$ is equivalent to $(t_1,\ldots,t_m,s_1,\ldots,\allowbreak
s_\mu)$, and hence to $(s_1,\ldots,s_\mu,1,\ldots,1)$. So any two
generating vectors of the same length become equivalent after $\mu$
stabilizations.
\end{proof}

\begin{proposition}
If $m>\ell(G)$, then any two generating $m$-vectors for $G$ are equivalent.
\label{prop:supraminimal generating sets}
\end{proposition}

\begin{proof}
Fix a generating vector $S=(s_1,\ldots,s_\mu)$ of minimal length, and let
$T=(t_1,\ldots,t_m)$ be any generating vector with $m$ elements. We will
show that $T$ is equivalent to $(s_1,\ldots,s_\mu,1,\ldots,1)$.

Put $G_i=\langle t_1,\ldots, t_i\rangle$, so $G=G_m\supseteq
G_{m-1}\supseteq \cdots \supseteq G_1\supseteq \{1\}$. Since $m>\ell$,
$G_i=G_{i-1}$ for some $i>0$, so $t_i$ can be written as a word in
$t_1,\ldots\,$, $t_{i-1}$. Therefore $T$ is equivalent to
$(t_1,\ldots,t_{i-1},1,t_{i+1},t_{i+2},\ldots,t_m)$. So we may assume that
$t_1=1$. Since $t_1=1$, $T$ is equivalent to $(s_1,t_1,\ldots,t_m)$. Let
$G_1=\langle s_1\rangle$ and $G_i=\langle s_1,t_1,\ldots,t_{i-1}\rangle$
for $i>1$. Since $s_1\neq1$, we must have $G_i=G_{i-1}$ for some $i>1$. So
$t_{i-1}\in \langle s_1,t_1,\ldots,t_{i-2}\rangle$, and therefore $T$ is
equivalent to $(s_1,t_1,\ldots,t_{i-2},1,t_i,\ldots,t_m)$. We may reselect
notation so that $T$ is equivalent to $(s_1,1,t_3,\ldots,t_m)$.
Inductively, suppose that $T$ is of the form $\{s_1,\ldots,s_k,t_{k+1},
\ldots, t_m\}$, for some $k<\mu$. Put $G_i=\langle s_1,\ldots,s_i\rangle$
for $i\leq k+1$ and to $\langle s_1,\ldots,s_k,s_{k+1},t_{k+1},\ldots,
t_{i-1}\rangle$ for $i>k+1$. For $i\leq k+1$, since $S$ is a minimal
generating set, we have $G_i\neq G_{i-1}$, so we must have $G_j=G_{j-1}$
for some $j>k+1$. So $t_{j-1}\in\langle s_1,\ldots,s_{k+1},t_{k+1},
\ldots,t_{j-2}\rangle$. This implies that $T$ is equivalent to
$(s_1,\ldots,s_{k+1},t_{k+1}, \ldots, t_{j-2},1,t_j, \ldots,t_m)$, and,
perhaps after reselecting nota\-tion, to
$(s_1,\ldots,s_{k+1},t_{k+2},\allowbreak \ldots,\allowbreak t_m)$. So $T$
is equivalent to $(s_1,\ldots,s_\mu,t_{\mu+1},\ldots,t_m)$, and hence
to~$(s_1,\ldots,s_\mu,1,\ldots,1)$.
\end{proof}

Since $\ell(G)+1$ is at most $1+\log_2(|G|)$,
proposition~\ref{prop:supraminimal generating sets} improves theorem~3
of~\cite{Gilman}, which shows that any two generating vectors of length at
least $2\log_2(|G|)$ are equivalent.

\begin{corollary}
If $g\geq 1+|G|\,\ell(G)$, then any two actions of $G$ on a handlebody of
genus $g$ are equivalent.
\label{coro:limit stabilization}
\end{corollary}

If $A=(\Z/p)^k$, for $p$ prime, then $\mu(A)=\ell(A)=k$. Thus by
theorem~\ref{thm:abelian}, if $p\geq5$, there are inequivalent actions on
the handlebody of genus $1+(\ell(A)-1)|A|$, showing that the estimate in
corollary~\ref{coro:limit stabilization} is the best possible, in
general. On the other hand, it appears to be far from the best possible for
many cases. Frequently there is a large gap between $\mu(G)$ and $\ell(G)$
(for example, all symmetric groups can be generated by two elements, but
have values of $\ell(G)$ that are arbitrarily large).

\section[Questions]{Questions}\label{sec:questions}
Our results on free actions are far from complete. The most obvious
question is:

\medskip\noindent \emph{Question 1:} Are all actions on genera above the
minimum genus equivalent?\medskip

That is, is $e(k)=1$ for all $k\geq1$?  Algebraically, if $n>\mu(G)$ are
any two $n$-element generating sets of $G$ Nielsen equivalent?  According
to p.~92 of~\cite{Lyndon-Schupp}, this algebraic version was first asked by
F. Waldhausen.  It has been resolved negatively for infinite groups.  The
first example appears to be due to G. A. Noskov \cite{Noskov}, and general
constructions are given in \cite{Evans2}. An affirmative answer to
Question~1 for finite groups would imply affirmative answers to the next
two questions:

\medskip\noindent \emph{Question 2:} Is every action the stabilization of
a minimal genus action?\medskip

\noindent That is, is $\E(0)\to \E(k)$ surjective for all $k$?
Algebraically, is every generating vector $(s_1,\ldots,s_n)$
Nielsen equivalent to a generating vector of the form
$(t_1,\ldots,t_\mu,1,\ldots,1)$?\medskip

\medskip\noindent \emph{Question 3:} Are all actions of a group $G$ on a
handlebody of genus $g$ equivalent after an elementary
stabilization?\medskip

\noindent That is, is $\E(k)\to \E(k+1)$ always constant?  Algebraically,
are all generating vectors of the forms $(s_1,\ldots,s_n,1)$ and
$(t_1,\ldots,t_n,\allowbreak 1)$ equivalent?

We can ask whether our example in section~\ref{sec:solvable groups} is the
smallest of its kind.

\medskip\noindent \emph{Question 4:} Are there weakly inequivalent actions
of a nilpotent group on a handlebody of genus less than $8193$?

\bibliographystyle{amsplain}

\end{document}